\documentstyle[12pt]{article}

\title{Fast graphs for the random walker}

\author{B\'alint Vir\'ag}

\newtheorem{thm}{Theorem}
\newtheorem{pro}[thm]{Proposition}
\newtheorem{fact}[thm]{Fact}
\newtheorem{lem}[thm]{Lemma}
\newtheorem{cor}[thm]{Corollary}
\newtheorem{protorem}[thm]{Remark}

\newtheorem{protoeg}[thm]{Example}

\newenvironment{thm*}[1]{\par \trivlist
 \itemindent\parindent \item[\bf Theorem #1]
 \it\ignorespaces}{\endtrivlist}
\newcommand{\qed}{\hfill\mbox{$\framebox(5,5)[]{}$}}

\newenvironment{eg}{\begin{protoeg}\rm}{\end{protoeg}}
\newenvironment{rem}{\begin{protorem}\rm}{\end{protorem}}
\newenvironment{proof}{\par \trivlist
 \itemindent\parindent \item[\hskip\labelsep\sc Proof.]
 \ignorespaces}{\qed\endtrivlist}
\newenvironment{proofof}[1]{\par \trivlist
 \itemindent\parindent \item[\hskip\labelsep\sc Proof of #1.]
 \ignorespaces}{\qed\endtrivlist}

 \DeclareSymbolFont{AMSb}{U}{msb}{m}{n}
\DeclareSymbolFontAlphabet{\mathbb}{AMSb}
\newcommand\mnote[1]{}
\newcommand\be{\begin{equation}}
\newcommand\bel[1]{{\mnote{#1}}\be\label{#1}}
\newcommand\ee{\end{equation}}
\newcommand\lb[1]{\label{#1}\mnote{#1}}
\newcommand{\Gr}{{\mathcal G}}

\newcommand{\eps}{\varepsilon}

\newcommand{\comment}[1]{}
\newcommand{\walk}{$\{X_k\}$}

\newcommand{\R}{{\mathbb R}}
\newcommand{\Rnn}{{\R_{\ge 0}}}

\newcommand{\ev}{\mbox{\bf E}}
\newcommand{\pr}{\mbox{\bf P}}
\newcommand{\one}{{\mathbf 1}}
\newcommand{\as}{\mbox{\hspace{.3cm} a.s.}}
\newcommand{\dist}{\mbox{\rm dist}}

\newcommand{\supp}{\mbox{\rm supp}}
\newcommand{\pid}{\pi_{-z}}
\newcommand{\sm}{{\raise0.3ex\hbox{$\scriptstyle \setminus$}}}
\newcommand{\Vd}{V\sm\{z\}}

\begin{document}
\maketitle
\begin{abstract}
Consider the time $T_{oz}$ when the random walk on a weighted
graph started at the vertex $o$ first hits the vertex set $z$.
We present lower bounds for $T_{oz}$ in terms of the volume of
$z$ and the graph distance between $o$ and $z$. The bounds are
for expected value and large deviations, and are asymptotically
sharp. We deduce rate of escape results for random walks on
infinite graphs of exponential or polynomial growth, and
resolve a conjecture of Benjamini and Peres.
\end{abstract}

 \footnotetext[1]
{{Research partially supported by the Lo\`eve Fellowship, the
Clay Mathematics Institute Liftoff Program, and NSF grant
\#DMS-9803597}.}
 \footnotetext[2]
{{\it AMS} 2000 {\it subject classifications. } Primary 60J15,
60F15; secondary 60J80.}

\section{Introduction}\lb{s intro}
A weighted graph $G=(V,w)$ is a set of vertices with a
symmetric nonnegative function $w$ on $V\times V$; the edges of
$G$ are given by the support of $w$. The goal of this paper is
to give a lower bound for the hitting times of sets for
reversible Markov chains, that is, random walks on weighted
graphs. At each step, the walk chooses a neighboring site at
random with odds given by the edge weights, and then moves
there. The weight $w_x$ of a vertex $x$ can be defined as the
sum of the weights over all incident edges. Define the weight
$w_z$ of a vertex set $z$ as the sum of the weights of the
vertices in the set. Let $T_{oz}$ denote the first time the
walk, started at $o$ at time $0$, visits the vertex or vertex
set $z$.

Consider the simple random walk on the nearest neighbor graph
of the integers from $o=0$ to $z=n$ and unit edge weights. As
it is easily computed,
 $$
 \ev T_{oz}=n^2.
 $$
Could the walk be faster if we assigned different edge weights?
Consider the biased simple random walk on the same graph with
odds of going left and right given by $1$ and $g>1$. This can
be realized as a random walk on a weighted graph with weight
$g^n$ on the $n$th edge. We have
 $$
 \ev T_{oz} \sim n(g+1)/(g-1).
 $$
The price we had to pay for higher speed is a higher weight on
vertex $z$: the parameter $w_z/w_o$ is constant $1$ for the
unbiased walk and $g^{n-1}$ for the biased walk. This raises
the question whether walks with fixed parameter $w_z/w_o$ can
be faster. Perhaps surprisingly, the first example, unbiased
simple random walk, is not the fastest. Theorem \ref{t py}
below implies that
 $$
  \inf_G \ev T_{oz} \sim {n^2 \over \log n},
 $$
where the infimum is taken over all weighted graphs (with
possibly more complicated structure) with vertices $o$, $z$ at
distance $n$ and $w_z/w_o=1$. In contrast, the second example
is asymptotically the fastest: by Corollary \ref{ldx}, the
infimum of $\ev T_{oz}$ taken over all weighted graphs with
vertices $o$, $z$ at distance $n$ and $w_z/w_o=g^{n-1}$ is
asymptotic to $n(g+1)/(g-1)$.

Let $r_{oz}$ denote the effective resistance between vertices
or vertex sets $o$, $z$ when the graph $G$  is thought of as an
electric network and conductances are given by the edge
weights. We are ready to state the main theorem.
 \begin{thm} \lb{t main}
Let $o$ be a vertex and $z$ be a vertex set in a weighted graph
with $\dist(o,z)\ge n+1$ for some integer $n\ge 0$. Let $m_g$,
$I_g(a)$ denote the mean and the large deviation rate function
of the hitting time $T'_{01}$ for biased simple random walk on
the integers with odds $1$ and $g$ of going left and right,
respectively.
 Then
 \begin{eqnarray*}
 \ev T_{o z } &\ge& m_gn + 1, \\
 \pr(T_{o z } \le an+1) &\le& e^{-I_g(a)n}.
 \end{eqnarray*}
Here $g$ may be taken to be either {\list{}{\itemsep 0in
\topsep 1ex}
 \item[(a)] the $g>1$ solution of \  $(g-1)^2g^{n-2}=2w_z/w_o$,
 or
 \item[(b)] $(w_z r_{oz})^{1/n}$.
 \endlist}
\noindent The classical formulas for $m_g$ and $I(a)$ are
 \begin{eqnarray}\label{mean}
 m_g&=&\ev T'_{01}=(g+1)/(g-1),\\
 e^{-I_g(a)}&=&\frac{g}{a+1}\
{{\bigg(\frac{g}{{a^2}-1}\bigg)}^{\frac{a-1}{2}}}\
 {{\Big(\frac{2a}{g+1}\Big)}^a}. \label{ld rate}
 \end{eqnarray}
 \end{thm}
 The implicit formula for the bound in part (a) can be made
 explicit.
 \begin{fact}\lb{f bd}
Theorem \ref{t main} still holds if $g$ is replaced by a
greater quantity. Set $\alpha:=n^2w_z/w_o \vee e$, then an
upper bound for $g$ in part (a) is given by
 \bel{gofa}
g':=\left[{5\alpha \over (\log \alpha)^2}\right]^{1 \over n-2}.
 \ee
 \end{fact}

In its applications Theorem \ref{t main} is related to the
bound of Varopoulos an Carne (1985), and in some cases, as in
the corollaries below, it yields sharper results. Part (b) is
related to the classical expression for commute time (see
formula (\ref{Chandra})) in the sense that it ties hitting
times and resistance. Also note that Theorem \ref{t main}
concerns large deviations, and therefore the bounds are more
precise than what follows from the Brownian (or Central Limit
Theorem) scaling limit; it compares random walks on graphs
directly to biased simple random walk. A version of the
expected value bound of part (b) was published in an earlier
paper, Vir\'ag (2000). Lee (1994ab) has a solution for the
optimization problem for expected value in the case of simple
path graphs and its continuous analogue. Large deviation
questions in random trees are studied by Dembo, Gantert, Peres
and Zeitouni (2001).

 The most important step in
proving Theorem \ref{t main} is a comparison of Laplace
transforms.
 \begin{pro}\lb{p main}
Using the notation of Theorem \ref{t main} (a), (b),
respectively, the Laplace transform of $T'_{0n}+1$ dominates
the Laplace transform of $T_{o z }$, that is, for $\lambda \ge
0$ we have
$$\ev e^{-\lambda (T'_{0n}+1)} \ge \ev e^{-\lambda
T_{o z }}.$$
 \end{pro}

It is possible to take limits of Theorem \ref{t main} in many
directions of its two parameters. The following asymptotic
result shows that in expected value and large deviations, the
fastest graphs are the ones corresponding to the asymmetric
simple random walks.

\begin{cor}[Large deviations for walks in graphs]\lb{ldx}
\ \\
Let $C,g>1$ and let $1<a<m_g$. Then
\begin{eqnarray*}
 \inf_G \ev T_{oz} &=& m_gn + O(1), \\
 \sup_G
 \pr[T_{oz} \le a n] &=& e^{-I_g(a)n+o(n)},
\end{eqnarray*}
where the $\inf$ and $\sup$ are over all  weighted graphs $G$
with vertex $o$ and vertex set $z$ satisfying $n \le
\dist(o,z)$ and $w_z/w_o \le C g^d$. The functions $O(1)$ and
$o(n)$ depend on $n,C,g,a$ only.
\end{cor}

Corollary \ref{ldx} allows us to prove a conjecture of
Benjamini and Peres (see Peres (1999)), originally stated for
unweighted trees. For an infinite graph $G$ with a fixed vertex
$o$, denote $w_n$ the total weight on edges at distance $n$.
The {\bf exponential upper growth} of $G$ is defined as
$\limsup w_n^{1/n}$. Let $|v|$ denote the graph distance
between vertices $v$ and $o$.

\begin{cor}[Exponential growth and lim sup speed] \lb{speed}
Let $G$ be an infinite weighted graph with exponential upper
growth $g_0$, and let $g=g_0 \vee 1$. Then the random walk
\walk\ on the graph satisfies
 $$
 \limsup {|X_k|\over k} \le {g-1 \over g+1} \ \as
 $$
\end{cor}
Note that equality holds for biased simple random walks, and,
perhaps surprisingly, even in certain {\it recurrent} graphs
(Example \ref{recurrent}).

In another scaling, we have

\begin{thm}[Large deviations in graphs of polynomial growth]
\lb{t py} \ \\ Let $0<c,d$ and $0<\alpha<2/(p+2)$. Then
 \begin{eqnarray*}
 \inf_G \ev T_{oz} &\sim& {2n^2 \over (p+2) \log n}, \\
  \sup_G \pr[T_{oz}<\alpha
n^2/\log n] &=& n^{-(\alpha(p+2)-2)^2/(8\alpha)+o(1)}.
 \end{eqnarray*}
The $\inf$ and $\sup$ are taken over graphs $G$ with vertices
$o$, $z$ satisfying $n \le \dist(o,z)$ and $w_z/w_o<c n^p$. The
function $o(1)$ converges to $0$ as $n\rightarrow \infty$ and
depends on $n$, $\alpha$, $c$, $p$ only.
\end{thm}

Surprisingly, Theorem \ref{t main} can be used to get sharp
results in this polynomial scaling, which is unlike the usual
domain for large deviation type bounds. Theorem \ref{t py}
implies a version of Khinchin's Law of the Iterated Logarithm
for random walks on infinite graphs. We say an infinite graph
$G$ has {\bf polynomial boundary growth} with power $p$ if $w_n
\le C n^p$ for all $n$ and fixed $C$.

\begin{cor}[Law of the single logarithm for walks on graphs]
 \lb{c ll1}
 \ \\
 We have
 $$
 \sup_G \left(\limsup {|X_k|\over \sqrt{k\log k}} \right)
 = {\sqrt{p+2}\over 2},
 $$
where the supremum is taken over random walks $\{X_k\}$ on
weighted graphs $G$ of polynomial boundary growth with power
$p$. The $\limsup$ is taken in the almost sure sense.
 \end{cor}
from known examples and the bounds of Varopoulos and Carne (1985).
Barlow and Perkins (1989) construct an unweighted subtree of 
${\bf Z}^2$ 
where the rate of escape is, up to a constant, the same as in 
Corollary \ref{c ll1}. 

In Sections \ref{s outline}, \ref{s flow}, \ref{s dcomp},
\ref{s array} and \ref{s Laplace} we present a proof of
Proposition \ref{p main}. These sections contain the main ideas
of the paper, which are outlined in Section \ref{s outline}. In
Section \ref{s Laplace} it is showed that (b) of Proposition
\ref{p main} implies (b) of Theorem \ref{t main}. A bit of
extra work is needed to prove part (a) Theorem \ref{t main},
and this is done in Section \ref{s main}. Implications of
Theorem \ref{t main} for graphs of exponential growth
(Corollary \ref{ldx} and Corollary \ref{speed}) are discussed
in Section \ref{s exp}. The polynomial case, including Theorem
\ref{t py} and Corollary \ref{c ll1}, is studied in Section
\ref{s poly}.

\section{Outline of the proof} \lb{s outline}

The proof of Theorem 1 follows easily from the Laplace
domination statement from Proposition \ref{p main}. The first
step in the proof of this proposition is to interpret the
Laplace transform probabilistically. But first, let us make
some conventions and introduce some notation.

For simplicity, the vertices in $z$ may be identified as a
single vertex (still denoted $z$), as it will not change any of
the quantities compared. A simple restriction argument also
shows that it suffices to prove the claim for finite weighted
graphs $G$. We assume further that all vertices in $V\sm \{z\}$
are accessible from $o$ without passing through $z$.

Let us call the object of our study a {\bf stopped random walk
law, SRWL}, defined as a quadruple $(K,o,V,z)$, where $V$ is a
finite set of vertices, $K$ is a reversible transition kernel
on $V$, $o,\;z \in V$ are vertices at which the random walk
will be started and stopped, respectively. Since we are only
interested in the walk before it reaches the vertex $z$, let us
introduce the notation $K_z$ for the transition probabilities
of the walk killed at $z$, that is
$K_z(x,y):=K(x,y)\one(x\not=z)$.

The Laplace transform of  $T_{oz}$ has the following
probabilistic interpretation. Let $\beta \in (0,1]$, and
consider the random walk which moves as the walk defined by $K_
z $, but is killed before each step with probability $1-\beta$.
Denote its kernel $K_\beta(x,y):=\beta K_ z (x,y)$. Let
$S_\beta$ denote the probability that this walk survives to hit
$z$. Then
 $$
S_\beta= \sum_{k=0}^\infty \pr(T_{o z }=k) \beta^k = \ev
\beta^{T_{o z }}.
 $$
This means that $S_\beta$, as a function of $-\log \beta$, is
the Laplace transform of $T_{o z }$.

The main difficulty in the proof of the proposition is that one
has to do optimization over a complicated geometric structure,
a graph. There are, however, graphs for which the optimization
is fairly simple, for example for the SRWL {\bf supported on a
path}. The definition is that the corresponding graph structure
is a finite path with $o,\;z$ at the two endpoints. Another
example is that of a {\bf dead-end RW law}, defined as a
stopped random walk law for which the probability of getting to
$ z $ is $0$. For such walks the optimization is trivial, since
the parameter $S_\beta$ is identically $0$.

Fortunately, every SRWL  can be ``decomposed'' into such simple
SRWLs in way that is convenient for our problem. First we
define some parameters that are natural for such decomposition.
Let $R_\beta$ denote the number of times the it visits $o$. Let
$\Gamma_\beta=R_\beta w_z/w_o$. The parameter $w_zr_{oz}$ of
part (a) will not enter directly into our analysis, but through
the bound
 \bel{Gamma wr}
 \Gamma_\beta = R_\beta w_z/w_o \le R_1 w_z/w_o =w_z  r_{o z }.
 \ee
The inequality here is trivial, where the second equality comes
from the well-known connection between random walks and
electric networks. Indeed, the probability that a random walk
on a graph started at vertex $o$ visits vertex $z$ before
returning to $o$ is given by $1/(w_o r_{oz})$. The number of
hits to $o$ before hitting $z$ is therefore a geometric random
variable with success probability $1/(w_o r_{oz})$, so its
expected value is $R_1=r_{oz} w_o$, just what we needed.

In short, Proposition \ref{p main} amounts to a comparison of
the parameter $S_\beta$ with the parameter
$\Gamma_\beta/R_\beta$ (part (a)), and with the parameter
$\Gamma_\beta$ (part (b)). These parameters are natural because
of
\def\mas{{\mathcal K}}
\begin{pro}[Decomposition of SRWLs]
\lb{dcomp} \ \\ Let $\mas$ be a SRWL and let $\beta \in (0,1)$.
There exist SRWLs $\{\mas_i\}_{i\in \Pi \cup \{o\}}$ and a
probability distribution $\alpha$ on $\Pi \cup \{o\}$ so that
$\{\mas_\pi\}_{\pi\in \Pi}$ are supported on paths, $\mas_o$ is
a dead-end RW law,
 {\list{}{\itemsep
0in \topsep 1ex}
\item[-]
the parameters $\Gamma_\beta$, $S_\beta$, and $R_\beta$ of
$\mas$ equal the convex combination with coefficients
$\alpha_i$ for the corresponding parameters of the $\mas_i$,
and
\item[-]
$\dist(o, z)$ in $\mas$ is not more than the corresponding
distance in the $\mas_i$.
\endlist}
\end{pro}
This proposition essentially says that the optimization can be
done on convex combinations of SRWLs on trivial graphs. As it
turns out, even this case is not completely straightforward,
especially for part (a). The analysis is presented in Sections
\ref{s array}, \ref{s Laplace} and \ref{s main}.

The proof of Proposition \ref{dcomp} depends on a duality
between stopped random walk laws and certain loss flows
presented in the next section.

\section{Random walks and flows}\lb{s flow}

 Proposition \ref{dcomp} claims that all SRWLs can be
replaced by convex combinations of basic SRWLs; this suggests a
representation of SRWLs as an elements of a linear space. It
will be a space of loss flows.

 Let $(K,o,V,z)$ be a SRWL, and let $\beta \in
(0,1)$. Consider the Green kernel $\Gr_\beta(x,y)$ defined by
$K_\beta$, which gives the expected number of times the walk
started at $x$ visits $y$: $\Gr_\beta(x,y):=\sum_{n=0}^{\infty}
K_\beta^n(x,y)$. Define the function $f:V\times V \rightarrow
\Rnn$
 \bel{d f}
 f(x,y):=\Gr_\beta(o,x)K_\beta(x,y),
 \ee
in words, the expected number of steps the walk defined by
$K_\beta$ makes from $x$ to $y$. This function encodes the
original random walk in a nice way. For example, it is easy to
see that reversibility of $K$ is reflected by  the fact that
for each cycle $\pi=(x_0,...,x_\ell=x_0)$ and its reversal
$\pi'$ the function $f$ satisfies
 \bel{reversible}
 f(\pi)=f(\pi').
 \ee
Here, and in the sequel, a {\bf function from $V\times V$
applied to a path} will mean the product of the values of the
function over the edges of the path.

Also, by comparison of the expected number of steps entering
and leaving a vertex $x\in V$ the following {\bf node law}
holds:
 \bel{Kirkhoff}
 \beta(f(V,x)+\one(o=x))\one(x\not= z ) = f(x,V).
 \ee
We refer to $f$ as a ``loss flow'' because it satisfies
Kirkhoff's node law for flows except for the factor $\beta<1$.

It is also possible to reconstruct the random walk from the
loss flow. Given a vertex set $V$, vertices $o,\;z$, a real
$\beta\in(0,1)$ and a nonnegative function $f$ on $V^2$
satisfying (\ref{reversible}) and (\ref{Kirkhoff}), define the
transition kernel
 $$
 K_ z (x,y):=\frac{f(x,y)}{\beta(f(V,x)+\one(o=x))}.
 $$
This kernel corresponds to a SRWL $(K,o,V,z)$ for which the
flow defined by (\ref{d f}) is $f$.

The relevant parameters $S_\beta$, $R_\beta$, $\Gamma_\beta$
can be expressed using the function $f$. Clearly:
 \bel{S 2}
 S_\beta = \sum_{x\in V}f(x, z ), \ \ \ \
 R_\beta = 1+\sum_{x\in V}f(x, o ).
 \ee
The parameter $\Gamma_\beta$ is somewhat harder to express.
Notice that for $x\in \Vd$, $y\in V$ the definition (\ref{d f})
and the fact that $K_\beta(x,y)=\beta w(x,y)/w_x$ implies $$
f(x,y)=\Gr_\beta(o,x)\beta w(x,y)/w_x, $$ and therefore for
$(y,x)\in \supp(f)$
 \bel{theta}
 \theta(x,y):=\frac{f(x,y)}{f(y,x)}
 \ee
 satisfies
 \bel{theta 2}
 \theta(x,y)=\frac{\Gr_\beta(o,x)w_y}
{\Gr_\beta(o,y)w_x}.
 \ee
Since $w(x, z )=w_x K(x, z )$, and $w_z=\sum_x w(x,z)$, the
parameter $\Gamma_\beta = \Gr_\beta(o,o) w_z/w_o$ can be
written as
 \bel{Gamma 2}
 \Gamma_\beta=\sum_x \Gr_\beta(o,o)w_x K(x, z )/w_o,
 \ee
where the sum runs over all neighbors of $z$. For every such
vertex $x$ we pick a simple path $\pi_x=(x_0=o,\ldots,
x_{\ell(x)}=x)$ for which $K(\pi_x)>0$. Repeated use of
equation (\ref{theta 2}) then transforms (\ref{Gamma 2}) to an
expression purely in terms of the flow $f$:
 \bel{Gamma 3}
 \Gamma_\beta
 =\sum_x \theta(\pi_x) \Gr_\beta(o,x)K(x, z )
 =\sum_x \theta(\pi_x) f(x, z )/\beta.
 \ee
We have expressed the three important parameters as functions
of the loss flow $f$.

\section{Decomposition of flows}

\lb{s dcomp} Consider a SRWL $(K,o,V,z)$, and let $f$ be a flow
defined by this SRWL (\ref{d f}). Fix $\theta$ as in
(\ref{theta}). Consider the set $F$ of nonnegative functions
$f_*$ on $V^2$ satisfying (\ref{Kirkhoff}, \ref{theta}) (with
$f$ replaced by $f_*$ in both), and
 $
 \supp(f_*)\subset \supp(f).
 $
 Note that (\ref{theta}) implies
(\ref{reversible}), so for each $f_*$ it is possible to define
a SRWL for which $f_*$ is the corresponding loss flow. Since
the parameters $\Gamma_\beta$, $S_\beta$, $R_\beta$ (\ref{Gamma
3}, \ref{S 2}) are clearly linear on $F$, the following lemma
will suffice for the proof of Proposition \ref{dcomp}.

\begin{lem}[Decomposition of loss flows]\lb{dec flow}\ \\
Let $\Pi$ be the set of simple paths from $o$ to $ z $. There
exists a probability distribution $\alpha$ on $\Pi \cup \{o\}$
so that
$$
f=\alpha_o f_o + \sum_{\pi\in \Pi} \alpha_\pi f_\pi,$$ where
$f_\pi \in F$ is supported on $\pi$ and $f_0\in F$ is supported
on $(\Vd)^2$.
\end{lem}

For the proof of this lemma, it is useful to know which $\pi\in
\Pi$ supports elements of $F$.
\begin{lem}\lb{path}
Let $\pi=(x_0=o,x_1,\ldots,x_\ell= z )$ be a simple path. There
exists $f_\pi \in F$ supported on $\pi$ if and only if
$\theta_i:=\theta(x_{i},x_{i-1})<\beta$ for all $1\le i\le
\ell$.
\end{lem}

\begin{proof}
There is a unique solution $f_\pi$ supported on $\pi$ for the
equations (\ref{Kirkhoff}) and (\ref{theta}). It can be
obtained inductively; set $\theta_0:=0$, then
 \begin{eqnarray*}
 f_\pi(x_{i-1},x_{i})&=&\prod_{j=1}^i
 \frac{\beta-\theta_{j-1}}{1-\beta\theta_{j}},\\
 f_\pi(x_{i},x_{i-1})&=&\theta_i f_\pi(x_{i-1},x_{i}).
 \end{eqnarray*}
The solution is nonnegative (equivalently, $f_\pi \in F$) if
and only if $\theta_i:=\theta(x_{i},x_{i-1})<\beta$ for all
$1\le i\le \ell$.
\end{proof}

\begin{proofof}{Lemma \ref{dec flow}}
$F$ is a closed, bounded, convex subset of a finite dimensional
vector space, so it equals the closed convex hull of its
extreme points (e.g. for a bit of overkill, by the Krein-Milman
Theorem). Thus it suffices to prove that all extreme points of
$F$ are supported on $(\Vd)^2$ or on simple paths. Indeed, let
$f_*$ be extreme point. Consider the directed graph on $V$
where $(x,y)$ is an edge iff $f_*(y,x)<\beta f_*(x,y)$
(equivalently, if $\theta(y,x)<\beta$ and $f_*(x,y)\not=0$).
Consider the set $V'$ of vertices which are connected to $ z $
by a path directed towards $ z $ in this graph.

First suppose that $o\notin V'$. Summing the node law
(\ref{Kirkhoff}) over elements of $V'$ yields
 \bel{inout1}
 \beta f_*(V,V') - \beta f_*(V, z ) = f_*(V',V).
 \ee
The definition of $V'$ implies that for $(x,y)\in (V\sm
V')\times V'$ we have $f_*(y,x)\ge \beta f_*(x,y)$; summation
yields
 \bel{inout2}
 f_*(V',V\sm V') \ge \beta f_*(V\sm V',V').
 \ee
Adding (\ref{inout2}), the trivial inequality $f_*(V',V') \ge
\beta f_*(V',V')$, and (\ref{inout1}) yields $0 \ge \beta
f_*(V, z )$ and therefore $\supp(f_*)\subset (\Vd)^2$.

The second case is when $o\in V'$, so there is a directed path
$\pi$ satisfying the assumptions of Lemma \ref{path} and that
$f_*(x,y)>0$ if $y$ follows $x$ in $\pi$. Thus there exists an
$f_\pi\in F$ supported on $\pi$, and for a small $\eps>0$, the
function $(1+\eps)f_*-\eps f_\pi$ is nonnegative hence an
element of $F$. As $f_*$ is an extreme point, $f_*=f_\pi$, and
the proof is complete.
\end{proofof}

\section{An array encoding a random walk law}\lb{s array}

The goal of this section is to extract the information in the
elementary SRWLs of the decomposition in Proposition
\ref{dcomp} into an array of numbers.

First assume that the chain is supported on a simple path
$\pi=(x_0=o,x_1,\ldots,x_\ell= z )$. Define the quantity
 \bel{d s}
 s(x,y):=\frac{\beta f(x,y)-f(y,x)}{f(x,y)-\beta f(y,x)}.
 \ee
We will express the relevant parameters in terms of the
$s(x,y)$. For $1\le i\le n$, the definition of $f$ implies that
$ s(x_{i-1},x_i) \in (0,\beta]$, and the node law
(\ref{Kirkhoff}) applied to $x_i$ implies that
 $$
 \beta f(x_{i-1},x_i) - f(x_i,x_{i-1}) =
  f(x_{i},x_{i+1})- \beta f(x_{i+1}, x_{i}) .
 $$
This makes the following product telescope:
 \bel{ftoz0}
  s(\pi) = \frac{\beta
  f(x_{\ell-1}, z )-f( z ,x_{\ell-1})}{f(o,x_1)-\beta f(x_1,o)}.
 \ee
 Since $f( z ,x_{\ell-1})=0$, we get $ s(x_{\ell-1}, z )=\beta$.
The node law (\ref{Kirkhoff}) applied to $o$ yields
 \bel{R s0}
 \beta(f(x_1,o)+1)=f(o,x_1),
 \ee
and this implies that the denominator of the right hand side of
(\ref{ftoz0}) also equals $\beta$. So if $\pid$ denotes the
path $\pi$ with its last vertex removed, then
 \bel{S s}
 S_\beta=f(x_{\ell-1}, z )= \beta  s(\pid).
 \ee
Note that the expression (\ref{R s0}) equals $\beta R_\beta$,
which, together with the definition of $s(o,x_1)$ yields
 \bel{R s}
 R_\beta = {1-\beta s(o,x) \over 1 - \beta^2}.
 \ee
Finally, we substitute (\ref{S s}) to (\ref{Gamma 3}) to get
 \bel{Gamma s}
 \Gamma_\beta = \theta(\pid) s(\pid)/\beta
 = h(s(\pid)) /\beta
 \ee
where
 \bel{d h}
 h(s):= s\frac{1- s\beta}{\beta- s}
 \ee
so that $h(s(x,y))=s(x,y)\times \theta(x,y)$. Now we turn to
the general case.

\begin{pro}[Array representation of SRWLs] \lb{p array} \ \\
Consider a SRWL, and let $\beta \in (0,1)$. There exists
{\list{}{\itemsep 0in \topsep 1ex}
\item[-]
a finite index set $\Pi$,
\item[-]
positive numbers $\alpha_\pi$ for $\pi\in\Pi$ with total sum at
most 1,
\item[-]
positive integers $\ell_\pi$ for $\pi\in\Pi$, with $\ell_\pi\ge
\dist(o,z)$ and
\item[-]
$s_{\pi,i}\in(0,\beta)$ for $\pi\in\Pi$, $1\le i< \ell_\pi$
\endlist}
so that the parameters of the SRWL satisfy
 \begin{eqnarray}
 S_\beta&=&
 \beta\sum_{\pi\in\Pi} \alpha_\pi
 \prod_{i=1}^{\ell_\pi-1} s_{\pi,i}, \label{S a}\\
 \Gamma_\beta&=&\sum_{\pi\in\Pi}
  \alpha_\pi \prod_{i=1}^{\ell_\pi-1}
  h(s_{\pi,i}),\label{Gamma a}\\
 R_\beta&\le&{2 \over 1-\beta^2}
  \left(1-\beta\sum_{\pi\in\Pi}
  \alpha_\pi s_{\pi_,1}\right). \label{R a}
 \end{eqnarray}
\end{pro}

\begin{proof}
We use the notation and the results of the decomposition in
Proposition \ref{dcomp}, so we can assume that our SRWL is a
convex combination of SRWLs supported on simple paths and a
dead-end RW law. For every simple path component
$\pi=(x^\pi_0,\ldots, x^\pi_{\ell_\pi})$, consider the stopped
random walk law there, and the flow defined in Section \ref{s
flow} for this walk. Let $s_{\pi,i}:=s(x^\pi_{i-1},x^\pi_{i})$,
as defined above (\ref{d s}).

The parameters $S\beta$ and $\Gr_\beta$ equal zero for the
dead-end random walk law component. Thus  (\ref{S a}) and
(\ref{Gamma a}) follow from linearity and the simple path case
(formulas (\ref{S s}, \ref{R s})).

The rest of the proof concerns the bound (\ref{R a}) for the
parameter $R_\beta$; it is relevant for part (a) Theorem \ref{t
main}, but not for part (b), and should be omitted at first
reading. For the dead-end random walk law component, the
$R_\beta$ is bounded above by the expected lifetime of the
walker $(1-\beta)^{-1}$ (in fact, this is sharp, achieved when
the graph consists of the vertex $o$ and a self-loop; if we
outlaw self-loops, the sharp bound becomes $(1-\beta)^{-2}$).
From this and the simple path case (\ref{R s}) we get
 $$
 R_\beta \le {\alpha_o\over 1-\beta} + \sum_{\pi \in \Pi}
 \alpha_\pi{1-\beta s_{\pi,1}  \over 1-\beta^2}.
 $$
Unfortunately, because of the possible self loops, this
expression for $R_\beta$ is messy, making the solution of the
optimization problem messy, too. To avoid this, we sacrifice
sharpness for simplicity, bounding the $(1-\beta)^{-1}$ and
$(1-\beta^2)^{-1}$ terms by $2(1-\beta^2)^{-1}$. This yields
the bound (\ref{R a}).
\end{proof}

\section{Laplace domination} \lb{s Laplace}

In this section we complete the proof of Proposition \ref{p
main} part (b) and Theorem \ref{t main} part (b) outlined in
Section \ref{s outline}. We will use the notation and results
introduced above.

\begin{proofof}{Proposition \ref{p main}, Part (b)}
We have seen in Section \ref{s outline} that it suffices to
bound $S_\beta$ in terms of $\Gamma_\beta$. Using the notation
and results of Proposition \ref{p array}, we can write
 \begin{eqnarray*}
  \Gamma_\beta
 &=&\sum_\pi
 \alpha_\pi \prod_{i=1}^{\ell_\pi-1} h( s_{\pi,i})\\
 &\ge&
 \sum_\pi
 \alpha_\pi h
  \left[\left(\prod_{i=1}^{\ell_\pi-1}
  s_{\pi,i}\right)^{\frac{1}{\ell_\pi-1}}
  \right]^{\ell_\pi-1}\\
 &\ge&
 \sum_\pi \alpha_\pi h\left[
 \left(\prod_{i=1}^{\ell_\pi-1}
 s_{\pi,i}\right)^{1/n} \right]^n\\
 &\ge&
 h\left[\left(\sum_\pi \alpha_\pi \prod_{i=1}^{\ell_\pi-1}
 s_{\pi,i}\right)^{1/n} \right]^n =
 h\left[\left(S_\beta/\beta\right)^{1/n}\right]^n.
 \end{eqnarray*}
The first inequality follows from Jensen's inequality and the
fact that $y\mapsto \log(h(e^y))$ is convex for $y\le \log
\beta$. The second, from the fact that the function $y\mapsto
h(y^{1/n})^n$ is increasing in $n$ for $y\in [0,1]$, and that
$\ell_\pi\ge n+1$. The third inequality follows from Jensen's
inequality and the fact that $y\mapsto h(y^{1/n})^n$ is convex
in $y$ for $y>0$.

Solving the above inequality for $S_\beta$, and using the fact
(\ref{Gamma wr}) that $g=(w_{ z }r_{o z })^{1/n}\ge
\Gamma_\beta^{1/n}$
 we get
 \bel{Laplace}
 S_\beta \le \beta
 \left( \frac{g+1-\sqrt{(g+1)^2-4\beta^2g}}{2\beta} \right)^n=
 \ev \beta^{T'_{0n}+1}.
 \ee
Note that $T'_{0n}$ is the sum of $n$ independent copies of
$T'_{01}$. Conditioning on the first step yields
 $
 \ev \beta^{T'_{01}}=\beta(g + (\ev \beta^{T'_{01}})^2)/(g+1)
 $,
and solving this equation gives the equality in
(\ref{Laplace}), a standard result. Thus the inequality in
(\ref{Laplace}), in terms of $-\log \beta$, is a comparison of
the Laplace transforms of $T_{o z }$ and $T'_{0n}+1$, as
required.
 \end{proofof}

\begin{proofof}{Theorem \ref{t main} Part (b)}
The expected value inequality follows from differentiating the
Laplace transforms at 0. For the large deviation inequality,
note that
 $$
\pr(T_{o z }-1\le an)\le  e^{\lambda an}\ev e^{-\lambda (T_{o z
}-1)}
 $$
for every $\lambda>0$ by Markov's inequality. Replacing the
Laplace transform on the right by that of $T'_{0n}$ we get
\begin{eqnarray*}
\pr(T_{o z }-1\le an) &\le&
 \inf_{\lambda>0} \ev e^{-\lambda (T'_{0n}-an)}\\
&=& (\inf_{\lambda>0} \ev e^{-\lambda (T'_{01}-a)})^n =
e^{-I(a)n}. \ \ \ \ \ \
\end{eqnarray*}
For the last equality, we used the fact that the infimum over
$\lambda\in \R$ is achieved when $\lambda>0$; this can be
checked by direct calculation.

Direct computation, or the law of large numbers, implies the
expression $(\ref{mean})$ for $m_g$, and a standard computation
using the Laplace transform of $T'_{01}$ yields its large
deviation rate function $I_g$ given by (\ref{ld rate}).
 \end{proofof}

\section{Proof of the main theorem, part (a)}\lb{s main}

In this section we prove Theorem \ref{t main}, part (a). It
suffices to prove Proposition \ref{p main}, part (a), a
comparison of Laplace transforms.  Given that, the proof of
Theorem \ref{t main}, part (b) in the previous section also
implies part (a).
 The optimization needed here is much more complicated; the most
technical part is presented separately at the end of the
section in Lemma \ref{hard-part}.

\begin{proofof}{Proposition \ref{p main}, part (a)}
Using the notation and results of Proposition \ref{p array}, we
can write
 \begin{eqnarray}
  \Gamma_\beta
 &=&\sum_\pi \alpha_\pi
 \prod_{i=1}^{\ell_\pi-1} h( s_{\pi,i}) \nonumber
 \\
 &\ge&
 \sum_\pi \alpha_\pi
 h(s_{\pi,1}) h\left[\left(\prod_{i=2}^{\ell_\pi-1}
  s_{\pi,i}\right)^{\frac{1}{\ell_\pi-2}}
 \right]^{\ell_\pi-2} \nonumber
\\
 &\ge&
\sum_\pi \alpha_\pi h(s_{\pi,1}) h(s_{\pi,*})^{n-1} \label{b
Gamma},
 \end{eqnarray}
where
 $$
s_{\pi,*}:=\left( \prod_{i=2}^{\ell_\pi-1}
s_{\pi,i}\right)^{1\over {n-1}}.
 $$
The first inequality follows from Jensen's inequality and the
fact that $y\mapsto \log(h(e^y))$ is convex for $y\le \log
\beta$. The second, from the fact that the function $y\mapsto
h(y^{1/n})^n$ is increasing in $n$ for $y\in [0,1]$, and that
$\ell_\pi\ge n+1$.

We will keep the parameter $ S_\beta=\beta\sum_{\pi\in\Pi}
\alpha_\pi
  s_{\pi,1} s_{\pi,*}^{n-1}
$ fixed and try to minimize the lower bound (given by (\ref{b
Gamma}) and (\ref{R a}) of Proposition \ref{p array})
 \bel{b wzwo}
 {w_z \over w_o}={\Gamma_\beta \over R_\beta} >
\left({1-\beta^2\over 2}\right)\frac{\sum_\pi \alpha_\pi
h(s_{\pi,1}) h(s_{\pi,*})^{n-1}}
{1-\beta\sum_{\pi\in\Pi}\alpha_\pi s_{\pi_,1}}
 \ee
as the parameters $s_{\pi,1}$ and $s_{\pi,*}$ range over the
set $[0,\beta)$. We first claim that the infimum is achieved on
this set. If $s_{\pi,1}$ converges to $\beta$, then
$h(s_{\pi,1})$ will converge to $\infty$, so the lower bound in
(\ref{b wzwo}) can only converge to a small value if
$h(s_{\pi,*})$ converges to $0$, in which case $s_{\pi,1}=0$ is
a better choice. The same argument can be made with the roles
of $s_{\pi,1}$ and $s_{\pi,*}$ reversed, so the infimum must
indeed be achieved on this set.

Lemma \ref{hard-part} below, where the hard part of the
optimization is done, implies that $s_{\pi,1}$ (respectively,
$s_{\pi,*}$) have to be the same for every $\pi$, so we may
drop the indices $\pi$. The lower bound (\ref{b wzwo}) reduces
to
 \bel{b wzwo2} {w_z \over w_o} >
 \left({1-\beta^2\over 2}\right)
 \frac{\alpha h(s_1) h(s_*)^{n-1}}{1-\beta
\alpha s_1},
 \ee
and we have $S_\beta=\beta \alpha s_1 s_*^{n-1}$. Since
$h(s_1)/s_1$ is increasing in $s_1$, it is clear that
increasing $\alpha$ while keeping $\alpha s_1$ fixed will not
change $S_\beta$ but will decrease the numerator on the right
hand side of (\ref{b wzwo2}). Therefore the minimum is achieved
when $\alpha$ is maximal, so we may take $\alpha=1$. After
cancellations,  the bound (\ref{b wzwo2}) reduces to
 \bel{b wzwo4} {w_z \over w_o} >
 \left({1-\beta^2\over 2}\right)\frac{h(s_*)^{n-1}}{
 \beta/s_1-1}.
 \ee
We are left to minimize this while keeping $s_1s_*^{n-1}$
fixed. The solution is
 \begin{eqnarray}\nonumber
 s_1&=&
 \frac{1-\beta^2}{(s_*-\beta)^2+1-\beta^2}s_*,
 \\ \label{bests1s*}
 {1 \over \beta/s_1-1} &=& {1-\beta^2 \over
 (\beta/s_*-1)(1-\beta s_*)} =
 {1-\beta^2 \over (\beta/s_*-1)^2} h(s_*)^{-1}.
 \end{eqnarray}
Clearly, $s_1<s_*$, so we have $ S_\beta< \beta s_*^n $. If we
set $g=h(s_*)$ then this gives exactly the inequality in
formula (\ref{Laplace}). Therefore, to conclude the proof it
suffices to show that $g$ is bounded above by the $g_0>1$
solution of
 $2w_z/w_o=(g_0-1)^2g_0^{n-2}$. Equivalently,
it suffices to show that ${2w_z / w_o} \ge (g-1)^2 g^{n-2}$.
This follows if we combine the simple bound
 $$
 {1-\beta^2 \over \beta/s_* -1} > h(s_*)- 1 = g-1
 $$
with formulas (\ref{b wzwo4}) and (\ref{bests1s*}).
\end{proofof}

\begin{rem}
At the price of complicated and long computations, this proof
can be modified to get the exact graph that maximizes the
chance of survival $S_\beta$ with $w_z/w_o$ fixed for a given
parameter $\beta$. This graph depends on $\beta$, and is a
simple path graph, except that in some cases a self-loop
appears at the vertex $o$. There are three points at which we
sacrificed sharpness for simplicity: the bound (\ref{R a}) in
Proposition \ref{p array} (this essentially eliminated the need
for self-loops at $o$), the bound for $s_1< s_*$ in the proof
above, and the last inequality of the proof, which eliminated
the dependence on $\beta$. If we do not allow self-loops at
$o$, expected hitting time (the $\beta \rightarrow 1$ case) is
minimized in the graphs of Example \ref{e fast}.
\end{rem}

\begin{proofof}{Fact \ref{f bd}}
For $g>1$ the expression $(g-1)^2g^{n-2}$ is increasing in $g$,
and can be bounded below by $(\log g)^2g^{n-2}$, so it suffices
to prove that this expression is at least $2w_z/w_o$. We
substitute (\ref{gofa}):
 $$
(\log g)^2g^{n-2}= \left({1 \over n-2} \log\left[ {5\alpha
\over (\log \alpha)^2} \right]\right)^2 {5\alpha \over (\log
\alpha)^2},
 $$
and replace the last $5\alpha$ by $5n^2w_z/w_o$ to get the
lower bound
 $$
{w_z\over w_o}\left({n \over n-2} \log\left[ {5\alpha \over
(\log \alpha)^2} \right]\right)^2 {5 \over (\log \alpha)^2}. $$
This can be bounded below by $w_z/w_o$ times
 $$ 5\left[{\log
\alpha +\log 5 - 2\log \log \alpha \over \log \alpha}\right]^2,
 $$
which is easily checked to be at least $2$.
\end{proofof}

\begin{lem}\lb{hard-part}
Suppose that $x_1$, $y_1$, $x_2$, $y_2$ achieve the minimum of
the expression
 \bel{mi}
 \alpha_1 h(x_1)h(y_1)^m+\alpha_2 h(x_2)^mh(y_2)
 \ee
subject to the constraints
 \begin{eqnarray}
 \alpha_1 x_1y_1^m + \alpha_2 x_2 y_2^m  &=&  c_1,
 \label{co 1}
 \\
 \alpha_1 x_1 + \alpha_2 x_2 &\le& c_2,
 \label{co 2}
 \\
 0 \ \ \le\ \  x_1,x_2,y_1,y_2 &<&\beta,
 \nonumber
 \end{eqnarray}
where $h$ is defined in (\ref{d h}) and $\alpha_i$, $c_i$ are
positive constants. Then $x_1=x_2$ and $y_1=y_2$.
\end{lem}

\begin{proof}
 \ \\
 \indent {\it Step 1.} The minimum can only occur for $x_i\le
y_i$. Otherwise, we may define new values
$x_i'=y_i'=(x_iy_i^m)^{1/(m+1)}$, this will not violate the
constraints, and will decrease (\ref{mi}) by the convexity of
$x\mapsto\log h(e^x)$.

{\it Step 2.} Minimum must be achieved in the interior of
$[0,\beta)^4$. By step 1 and symmetry, the only other case is
$x_1=0$, for this we may assume $y_1=0$ (it makes no
difference). Let $c_2'$ be the value of the left hand side in
(\ref{co 2}). It is straightforward to check that the unique
solution of (\ref{co 1}) and
 \bel{co 2'}
 \alpha_1 x_1 + \alpha_2 x_2 = c_2'
 \ee
for which $y_1=y_2$ and $x_1=x_2$ gives a smaller value for
(\ref{mi}).

{\it Step 3.} Minimum is in fact achieved when equality holds
in (\ref{co 2}), but we do not need to show this. From now on
we will only use that (\ref{mi}) is also minimal when (\ref{co
2}) is replaced by the equality constraint (\ref{co 2'}), where
$c_2'$ is the actual value of the left hand side of (\ref{co
2}).

{\it Step 4.} Since we excluded the case that the minimum
occurs on the boundary, it can only occur where the 4
dimensional gradient of (\ref{mi}) is perpendicular to the 2
dimensional surface determined by (\ref{co 1}, \ref{co 2'}), or
the derivative evaluated at two linearly independent vectors
tangent to this surface is 0. This means that the Jacobian of
the map given by the function (\ref{mi}) and the left hand
sides of (\ref{co 1}, \ref{co 2'}) has a two-dimensional
nullspace, so any $3 \times 3$ submatrix must be singular. The
Jacobian is a $4\times 3$ matrix; the first and third columns
are computed as
 $$
 \left[
 \begin{array}{l}
 \alpha_1 h'(x_1)h(y_1)^m \\
 \alpha_1 y_1^m \\
 \alpha_1
 \end{array}
 \right],\ \
 \left[
 \begin{array}{llll}
 \alpha_1 mh(y_1)^{m-1}h'(y_1)h(x_1) \\
 \alpha_1 m x_1y_1^{m-1} \\
 0
 \end{array}
 \right].
 $$
We get the other two columns by replacing the index $1$ by $2$.
Now we set $r(x):=h(x)/x$ and divide the first two columns by
$\alpha_i$, and the last ones by entries in the second row:
 $$
 \left[
\begin{array}{llll}
h'(x_1)h(y_1)^m            &  h'(x_2)h(y_2)^m            &
r(y_1)^{m-1}h'(y_1)r(x_1)  &  r(y_2)^{m-1}h'(y_2)r(x_2) \\
y_1^m                      &  y_2^m                      &
1&1             \\
1&1&0&0
\end{array}
\right].
 $$
The determinant of the right $3\times 3$ submatrix has to be 0,
and this happens if and only if $f(x_1,y_1)=f(x_2,y_2)$ with
 $$
 f(x,y)=r(y)^{m-1}h'(y)r(x).
 $$
The left $3\times 3$ submatrix simplifies if we divide the
first row by $f(x_i,y_i)$:
 $$\left[
\begin{array}{lll}
y_1^m {h'(x_1)r(y_1)^m\over h'(y_1)h(x_1)} &
y_2^m {h'(x_2)r(y_2)^m\over h'(y_2)h(x_2)} & 1 \\
y_1^m                      &  y_2^m                      & 1 \\
1&1&0
 \end{array}
 \right].
 $$
After computing the determinant, we get that this matrix is
singular iff $g(x_1,y_1)=g(x_2,y_2)$ with
 $$
 g(x,y)=y^m\left(1-\frac{h'(x)}{r(x)}\frac{r(y)}{h'(y)} \right).
 $$
To complete the proof, we have to show that the map
$(x,y)\mapsto (f,g)$ is injective. This follows from the fact
that if $x\le y$, then $f$ is increasing in $x$, $y$ and $g$ is
decreasing in $x$ and increasing in $y$. Consider
$x_1,x_2,y_1,y_2$; the two interesting cases are $x_1<x_2,
y_1<y_2$, and $x_1<x_2,y_1>y_2$. In the first case
$f(x_1,y_1)<f(x_1,y_2)<f(x_2,y_2)$, in the second,
$g(x_1,y_1)>g(x_2,y_1)>g(x_2,y_2)$.
\end{proof}

\section{Results for graphs of exponential growth}\lb{s exp}

This section contains the proofs of Corollary \ref{ldx} and
Corollary \ref{speed}. We then give a recurrent example in
which the inequality of Corollary \ref{speed} is sharp.

Note that biased simple random walks achieve the bounds of
Corollary \ref{ldx} by the Law of Large Numbers and the Large
Deviation Principle. Thus it suffices to prove the following
Corollary to Theorem \ref{t main}. Its claim is more precise
than the lower bound of Corollary \ref{ldx}.

\begin{cor}\ Let $C,g>0$ and let $1<a_0<(g+1)/(g-1)$.
There exists $C_1,\;C_2$ so that
 \begin{eqnarray*}
 \ev T_{oz} &>& n(g+1)/(g-1) - C_2,
 \\
 \pr[T_{oz} \le a n] &<& C_1 e^{-I_g(a)n}
 \end{eqnarray*}
for any weighted graph with vertices $o$, $z$ satisfying $n \le
\dist(o,z)$ and $w_z/w_o \le C g^d$, and for all $a\in
[a_0,(g+1)/(g-1)]$.
\end{cor}

\begin{proof}
Set
\begin{eqnarray*}
n' &:=& n-1, \\
a' &:=& {an - 1 \over n-1}, \\
g' &:=&\left[{5(n-1)^2 C g^n\over [\log ((n-1)^2 C
g^n)]^2}\right]^{1/(n-3)}.
\end{eqnarray*}
Note that $g'>g$ as well as $a'>a$, and $a'$, $g'$ are bounded
by some constants $a'_{max}$ and $g'_{max}$ for all $n$.  Also,
for all large $n$, we have $1<a'<(g'+1)/(g'-1)$. For these $n$
we apply Theorem \ref{t main} and Fact \ref{f bd} with
parameters $n'$, $Cg^n$, $a'$ to get
 \begin{eqnarray}\label{q2 mn} \pr[T_{oz} \le an] &\le&
 e^{-n'I_{g'}(a')}, \\
 \ev T_{oz} &\ge& (g'+1)(g'-1)(n-1) + 1. \label{q3 mn}
 \end{eqnarray}
The definition of $g'$ implies that
 $
 {g'^{n-3} /  g^{n-3}} = O(1),
 $
 and therefore
 $
 g'-g = O(1/n).
 $
This and (\ref{q3 mn}) proves the proposed expected value
bound. Clearly $a'-a = O(1/n)$. The function $(g,a)\rightarrow
I_g(a)$ and its derivative are continuous on the set
$[g,g'_{max}]\times[a_0,a'_{max}]$, and hence bounded.
Therefore
 \begin{eqnarray}\nonumber
nI_g(a)-n'I_{g'}(a') &\le& n(I_g(a)-I_{g'}(a)) + c \\
 &\le& nc_1((g'-g)+(a'-a)) + c \le c_2. \label{derivative}
 \end{eqnarray}
The boundedness of $I_g(a)$ implies that we can ignore the
first few $n$ for the price of increasing the constant $C_1$.
Thus (\ref{q2 mn}) and (\ref{derivative}) imply the proposed
large deviation bound.
 \end{proof}

We are ready to prove Corollary \ref{speed}.
 \begin{proofof}{Corollary \ref{speed}}
Let $g'>g$ be arbitrary, and let $T_n$, $w_n$ denote the
hitting time and the total weight of the set of edges at
distance $n$ from $o$, respectively. Then
 \bel{limsup}
 \limsup |X_k|/k = \limsup n/T_n.
 \ee
 For all large $n$ we have $w_n/w_o<g'^n$, so by
Proposition \ref{s exp} each event $T_n\le an$ has probability
at most $c_1 e^{-c_2 n}$. By the Borel-Cantelli Lemma only
finitely many of these events happen. Thus (\ref{limsup}) is at
most $(g'-1)/(g'+1)$, and since $g'>g$ was arbitrary, the
Corollary follows.
\end{proofof}

\begin{eg}\lb{recurrent}
Let $g\ge 2$ be an integer, and consider the graph of the
nonnegative integers with $g$-ary trees of depth $d_i$ attached
at vertex $i$ for every $i$. If $d_i$ increases fast enough,
then by the time the walk started from $0$ visits the leaves of
the tree at $d_i$, its speed will be nearly as high as the
speed of the walk on the $g$-ary tree, and the upper growth of
this graph is just $g$. This gives an example of a recurrent
graph for which equality is achieved in Corollary \ref{speed}.
\end{eg}

\section{Hitting times in graphs of polynomial growth}\lb{s poly}

In this section we prove Theorem \ref{t py}. One direction of
the inequalities are simple corollaries to Theorem \ref{t
main}; the other direction is provided by the ``fast graphs''
of Example \ref{e fast} (expected value) and Example \ref{e
fast2} (large deviations).

In the end of the section, we prove Corollary \ref{c ll1}, a
version of Khinchin's Law of the Iterated Logarithm. Again, we
prove two inequalities; the first is provided by a Corollary to
Theorem \ref{t py}, the second, by Example \ref{e fast3}.

\begin{cor}
\lb{p pybd}
 \ \\
 Let $0<c,d$ and $0<\alpha<2/(p+2)$. We have
 \begin{eqnarray} \label{pybd ev}
 \ev T_{oz} &>&  {2n^2 \over (p+2)\log n} - O(n),
 \\
 \pr[T_{oz}\le\alpha n^2/\log n] &<&
 n^{-(\alpha(p+2)-2)^2/(8\alpha)+o(1)}
 \label{pybd}
 \end{eqnarray}
for all graphs $G$ with vertices $o$, $z$ satisfying $n \le
\dist(o,z)$ and $w_z/w_o\le c n^p$. The functions $o(1),\;O(n)$
depend on $n$, $\alpha$, $c$, $p$ only.
\end{cor}

 \begin{proof}
We apply Theorem \ref{t main}, part (a) to the graph in
question. The parameters we use are $n'=n-1$, $a_n$ which is
the solution of $a_nn'+1 = \alpha n^2/\log n$, and $g_n=
n^{(p+2)/n}$. This will cover the graphs in question, since
 $$
 (g_n-1)^2g_n^{n'} \sim
 \left({ 1\over n}\log n^{p+2}\right)^2n^{p+2}
 = (p+2) n^p \log n,
 $$
and this dominates $2w_z/w_o=2cn^p$ for large $n$. The theorem
yields
 $
 \ev T_{oz} >
 ( 1 + 2/(g_n-1))n'
 $,
and using the fact that $1/(g_n-1) = 1/\log g_n +O(1)$,  the
first claim (\ref{pybd ev}) follows.

Theorem \ref{t main} also yields the bound (\ref{ld rate}) on
$\pr[T_{oz}\le\alpha n^2/\log n]$  which we rewrite as follows:
 $$
 \left[1-{(g-1)^2\over (g+1)^2}\right]^{a \over 2}
 \cdot g^{1/2}
 \cdot \left[1+{1\over a^2-1}\right]^{a-1\over 2}
 \cdot \left[1-{1\over a+1} \right].
 $$
Substituting the parameters for our case and taking logarithms
we get
 $$
 {-a_n (\log g_n)^2 \over 8} + {\log g_n \over 2}
 + {a_n \over 2a_n^2} - {1\over a_n} + o(n^{-1}\log n).
 $$
Multiplying by $n$ and substituting the formulas for $a_n$ and
$\log g_n$ we get
 $$
\log n \left( {-\alpha(p+2)^2 \over 8} + {p+2 \over 2}  - {1
\over 2\alpha} +o(1)\right).
 $$
Exponentiation yields the bound (\ref{pybd}).
\end{proof}

\begin{eg} \lb{e fast}
We now show a family of fast simple path graphs of polynomial
growth. Let $g>1$, and consider the simple path graph with
vertices denoted $o=0,1,\ldots,n=z$, edges $e_i=(i-1,i)$ and
edge weights
\begin{eqnarray*}
 w(e_1) &=& 1, \\
 w(e_i) &=& (g-1)g^{i-2} \mbox{\ \ \ \ \ \ for }2\le i \le n-1,\\
 w(e_n) &=& (g-1)^2g^{n-3}.
\end{eqnarray*}
Heuristically, the walker has a positive drift when it is away
from the endpoints of the path; the price is large negative
drifts at the two ends.

Consider the stopped random walk on this graph and the flow
associated with parameter $\beta=1$ as defined in (\ref{d f}).
This flow is uniquely determined by the flow property
(\ref{Kirkhoff}) and the requirement that $f(i,i+1)/f(i,i-1) =
w(i,i+1)/w(i,i-1)$. These equations have solution:
$$\begin{array}{lcll}
 f(1,0)&=&g/(g-1)^2,& \\
 f(n,n-1)&=&0,& \\
 f(i,i-1)&=& 1/(g-1)& \hspace{4ex} \mbox{for }2\le i \le n-1,\\
 f(i-1,i)&=&f(i,i-1)+1& \hspace{4ex} \mbox{for }1\le i \le n.
\end{array}$$
And $f(i,j)=0$ elsewhere. Then clearly
 $$
 \ev T_{oz} = \sum_{i,j} f(i,j) = 2(n-2)/(g-1) + 2g/(g-1)^2 + n.
 $$
Now consider the case when
\begin{eqnarray} g &=& \left[{n^{p+2}
\over (\log (n^{p+2}))^2}\right]^{1/n},\mbox{\ \ \ \ \ so that}
\label{polyg}
\\ g-1 &\sim& \log g \ \sim \ {\log(n^{p+2}) \over n}. \nonumber
\end{eqnarray}
This example proves one direction of the expected value bound
in Theorem \ref{t py}, since
 \begin{eqnarray*}
  w_z/w_o &=& (g-1)^2 g^{n-3} \sim n^p,\\
  \ev T_{oz} &\sim& {2n^2 \over (p+2) \log n}.
 \end{eqnarray*}
\end{eg}

\begin{eg} \lb{e fast2}
We now show that in the previous example the large deviation
bounds of Proposition (\ref{p pybd}) are also achieved. It is
perhaps surprising that the bounds are sharp even in this
scaling. We want to estimate the probability that the hitting
time is short by dividing the path into three segments. Let
 $$
 m:=\lfloor n/\log n \rfloor, \ \ \ \ v:=n-1.
 $$
 We expect the walk to spend
most of its time between the vertices $m$ and $v$. More
precisely, let
 \begin{eqnarray*}
 t&:=&\alpha n^2/\log n, \\ t'&:=&t/(\log n)^{1/2}\ =\ \alpha
n^2 / (\log n)^{3/2},
 \end{eqnarray*}
 so by the strong Markov property
 \bel{smp}
 \pr[T_{oz}<t] \ \ge\
 \pr [T_{mv}<t-t']\; \pr[T_{om}+T_{vz}<t'].
 \ee
 The second factor can be bounded using the classical formula for
commute time (see Chandra et al (1989)) and Markov's
inequality. For any weighted graph and vertices $o,\;z$, if
$r_{oz}$ denotes effective resistance, then
 \bel{Chandra}
 \ev [T_{oz} + T_{zo}]= w_Vr_{oz}.
 \ee
By the series rule
\begin{eqnarray}
r_{om} &=& \sum_{i=1}^m w(e_i)^{-1}
  = 1 + {1-(1/g)^{m-1}
\over (g-1)(1-1/g)} \nonumber \\
 &\asymp& (g-1)^{-2} \ \asymp\  {n^2 /(\log n)^2}, \label{r0m}
\end{eqnarray}
the same way we get
 \bel{r0n}
 r_{oz} \asymp {n^2 /(\log n)^2},
 \ee
and $r_{vz} \sim n^{-p}$. Also, the total sum of edge weights
satisfies
 $$
 w_E = 1+ g^{n-3}(g-1)^2+(g-1)\sum_{i=2}^{n-1} g^{i-2}
 \asymp n^p + g^{n-2}
 \asymp {n^{p+2} \over (\log n)^2}.
 $$
 For the edges $E'$ on the path between vertices $o$, $m$ we have
 $$
 w_{E'} = 1 + (g-1)\sum_{i=2}^m g^{i-2}\asymp (g^{m-1}-1)\asymp
1,
 $$
so we can apply (\ref{Chandra}) twice:
 $$ \ev
[T_{om}+T_{vz}] \le 2w_{E'}r_{om} + 2w_{E}r_{vz} \asymp {n^2
\over (\log n)^2}.
 $$
Markov's inequality concludes the bound on the chance of the
complement of the last event of (\ref{smp}):
 $$
\pr[T_{om}+T_{vz} \ge t'] \le {\ev[ T_{om}+T_{vz}] \over t'} =
O(\log n)^{1/2}.
 $$
It remains to bound the first factor:
 $$
\pr[T_{mv}<t-t'] \ge \pr[T_{mv}<t-t'\ |\ T_{mv}<T_{mo}]
\pr[T_{mv}<T_{mo}].
 $$
The second term here can be computed using resistances:
 $$
\pr[T_{mv}<T_{mo}]={1 \over 1+r_{mv}/r_{mo}} \ge {1 \over
1+r_{ov}/r_{mo}} \ge c>0.
 $$
The constant lower bound follows from (\ref{r0m}) and
(\ref{r0n}). To bound the first term, first note that the walk
started at $m$ and conditioned on the event $T_{mv}<T_{mo}$ is
a Doob transform of the original walk, a reversible random walk
in which the forward drift is bounded below by the forward
drift in the original walk. Therefore, by stochastic
domination,
 $$
 \pr[T_{mv}<t-t'\ |\ T_{mv}<T_{mo}]\ \ge\ \pr[T'_{0,n-1-m}<t-t'],
 $$
where $T'$ denotes hitting time for biased simple random walk
$\{X'_k\}$ on the integers with odds of going left and right
equal $1:g$. The second probability is bounded below by the
probability of a smaller event, which in turn can be bounded
using Lemma \ref{tapp}:
 $$
 \pr[X'_{\lfloor t-t'-1 \rfloor}\ge n-1-m]  >
 n^{-(\alpha(p+2)-2)^2/(8\alpha)+o(1)}.
 $$
All together, this example gives one direction in the large
deviation bound of Theorem \ref{t py}:
 $$
 \pr[T_{oz} < \alpha n^2/\log n] >
 n^{-(\alpha(p+2)-2)^2/(8\alpha)+o(1)}.
 $$
\end{eg}

We now turn to the proof of the simple lemma we used in the
previous example. We were unable to locate a theorem in the
literature that would imply this claim.

\begin{lem}\lb{tapp}
Let $p \ge  0$  and $\{X_k\}$ be biased simple random walk on
the integers with odds for going left and right given by $1$
and $g=g(p,n)$ defined in formula (\ref{polyg}). Let
$\alpha<2/(p+2)$ and let
 $
 t=t(n) \sim \alpha n^2/\log n
 $.
Then
 $$
 \pr [X_t \ge n] > n^{-(\alpha(p+2)-2)^2/(8\alpha)+o(1)}.
 $$
\end{lem}

\begin{proof}
Without loss of generality we may assume that $t$ is even. Let
$n_1= n(1+1/ \log n)$. Then
 \begin{eqnarray}
 \pr[X_t \ge n] &\ge& \pr[X_t \in [n,n_1]]\nonumber \\
 &\ge& {n_1-n-2 \over 2}
 \inf_{n\le {m \atop \mbox{\tiny \ even\ }} \le n_1}
 \pr[X_t=m].
 \label{tail}
 \end{eqnarray}
We now use the binomial formula to get that for $t,\;m$ even
 $$ \pr[X_t=m] = { t
 \choose (t+m)/2 }{g^{(t+m)/2} \over (1+g)^t}.
 $$
 By Stirling's formula and the fact that
$t(n),\ m(n)\rightarrow \infty$ we get
 \begin{eqnarray*}
 2^{-t} { t \choose (t+m)/2 } &\sim&
 {1 \over \sqrt {\pi}}
 \frac{t^{t+1/2}}{(t+m)^{(t+m+1)/2}(t-m)^{(t-m+1)/2}}
 \\&\asymp&
 t^{-1/2}\left[1+ {m^2 \over t^2-m^2}\right]^{t/2}
 \left[1-{2m \over t+m} \right]^{m/2}.
 \end{eqnarray*}
 The remaining factor can be written as
 $$
 2^{t}{g^{(t+m)/2} \over (1+g)^t} =
 g^{m/2}\left[1 - {(1-g)^2 \over (1+g)^2}\right]^{t/2}.
 $$
 Using the expression (\ref{polyg}) for
$g$, the fact that $g(n) \rightarrow 1$ and that $m(n)/t(n)
\rightarrow 0$ we get that
 \begin{eqnarray*}
 \log \pr[X_t=m] &\sim&
 -{\log t \over 2}
 + {m^2 \over t^2}\cdot {t\over 2}
 - {2m \over t} \cdot{m\over 2}
 + {m \over 2} \log g
 - {(\log g)^2 \over 4}\cdot {t\over 2}
 \\ &\sim& \left(
 -1
 +1/(2\alpha)
 -1/\alpha
 + (p+2)/2
 - (p+2)^2\alpha/8
 \right) \log n
 \\  &=&
 \left(-1-(\alpha(p+2)-2)^2/(8\alpha)\right) \log n.
 \end{eqnarray*}
The convergence is uniform over all $m \in [n,n_1]$. This and
(\ref{tail}) imply the claim of the lemma.
\end{proof}

We now turn to the proof of the graph version of Khinchin's Law
of the Iterated Logarithm. The upper bound is a Corollary to
Theorem \ref{t py}.

\begin{cor}[Law of the single logarithm, upper bound] \lb{c ll2}
\ \\ For random walks $\{X_k\}$ on weighted graphs with
polynomial boundary growth with power $p$ we have
 \bel{e ll}
\limsup {|X_k|\over \sqrt{k\log k}}  \le {\sqrt{p+2}\over 2} \
\as
 \ee
 \end{cor}

\begin{proof}
Let $a=2/(p+2)$, let $a''<a'<a$, and let
 \bel{d f()}
 f(t)=\sqrt{(t \log t)/(2a'')}.
 \ee
 Let $m>1$ an integer, for every $k$, let $\ell_k$
denote the distance of the farthest vertex visited up to time
$k$, and let $\ell'_k$ be the greatest integer so that
$\ell'^m_k<\ell_k$. Then
 $$
{|X_k|\over f(k)} \le {\ell_k \over f(T_{\ell_k})} \le
{(\ell'_k+1)^m \over f(T_{\ell'^m_k})} = {(\ell'_k+1)^m \over
\ell'^m_k} {\ell'^m_k \over f(T_{\ell'^m_k})}.
 $$
Taking lim sup we get
 \bel{khin} \limsup_{k\rightarrow \infty}
{|X_k|\over f(k)} \le \lim_{k\rightarrow \infty} {(\ell'_k+1)^m
\over \ell'^m_k} \limsup_{k\rightarrow \infty} {\ell'^m_k \over
f(T_{\ell'^m_k})} =\limsup_{\ell\rightarrow \infty} {\ell^m
\over f(T_{\ell^m})}.
 \ee
We are taking $m$th powers to make a sequence of probabilities
summable. Now consider the function
 \bel{d g()}
 g(\ell)=a'\ell^2/\log \ell,
 \ee
an upper bound for  $f^{-1}$, so that $t<g(f(t))$ for all large
$t$. For all large $\ell$  $g(\ell)$ is increasing, and we have
 \begin{eqnarray*}
\pr[f(T_\ell)< \ell] &=& \pr[g(f(T_\ell))< g(\ell)]
\\ &\le& \pr[T_\ell < g(\ell)]
\\ &\le& \ell^{-(a'(p-2)-2)^2/(8a')+o(1)}.
 \end{eqnarray*}
The last inequality follows from Theorem \ref{t py}. For large
$\ell$ the right hand side is bounded above by $\ell^{-c}$ for
some $c>0$, so it is summable over the subsequence of
$m$-powers if $cm>1$. Therefore by the Borel-Cantelli lemma
$f(T_{\ell^m})\ge\ell^m$ eventually a.s., so the right hand
side of (\ref{khin}) is at most 1. Since $a''<a$ was arbitrary,
the corollary follows.
 \end{proof}

\begin{eg} \lb{e fast3}
Using Example \ref{e fast} it is easy to construct an example
for the sharpness of Corollary \ref{c ll2}, and thus prove
Corollary \ref{c ll1}. Let $x_i$ be a sequence where
$x_i-x_{i-1}$ is positive and rapidly increasing. Consider the
sequence of simple path graphs $G_i$ of length
$n_i=x_i-x_{i-1}$ and of polynomial growth $w_z/w_o=n_i^p$
constructed in Example \ref{e fast}. We concatenate them in
increasing order to get an infinite simple path graph. By
picking $x_i-x_{i-1}$ to be rapidly increasing, it can be
achieved that the dominant term in the expected hitting time
$\ev T_{0,x_i}$ will be the expected hitting time $\ev T_{oz}$
in the graph $G_i$. This means that if we set $a=2/(p+2)$ the
walk in the concatenated  graph has
 $$
 \ev T_{0,x_i} \sim a{x_i^2 \over \log x_i}.
 $$
Let $a''>a'>a$. With $g$ as in (\ref{d g()}), by Markov's
inequality, for all large $i$ and all $y<x_i$ we have
 $$ \pr\left[T_{y,x_i} \le g(x_i)\right]
 \ge 1-a/a'.
 $$
The function $f$ (\ref{d f()}) is a lower bound for the inverse
of $g$ in the sense that $f(g(\ell))<\ell$ for all large
$\ell$. Therefore for large $i$
 $$
 \pr\left[f(T_{y,x_i}) \le x_i\right] \ge 1-a/a'.
 $$
The following implications are simple:
 $$\limsup_{i \rightarrow \infty}
 \left\{{x_i  \over f(T_{x_i})} \ge 1
 \right\}
 \Rightarrow
 \limsup_{k \rightarrow \infty}
 \left\{ {X_k \over f(k)}  \ge 1 \right\}
 \Rightarrow
 \left\{\limsup_{k \rightarrow \infty}
 { X_k \over f(k)}  \ge 1 \right\}.
 $$
The first event has probability at least $1-a'/a$ no matter
which vertex the walk is started at. Let $A$ denote the last
event; we then have $\pr A \ge 1-a'/a$ even if we start the
walk at a different time (as opposed to time 0). Thus by
L\'evy's 0-1 law
 $$
 1-a'/a\le\pr[A|X_0,\ldots, X_k]
 \rightarrow \one_A  \ \ \ a.s.
 $$
Thus $\pr A=1$, and since $a''>a$ was arbitrary, the lower
bound in Corollary \ref{c ll1} follows.
\end{eg}

 \noindent {\bf Acknowledgments.}
The author thanks Noam Berger, Russell Lyons and Yuval Peres for helpful
comments on previous versions.

\section*{References}
\newcommand{\andd}{and\ }
\def\labelenumi{[\arabic{enumi}]}
\begin{enumerate}
 \item
Barlow, M. T., Perkins, E. A. (1989) Symmetric Markov chains in 
${\bf Z}^d$: how fast can they
move? {\it Probab. Theory Relat. Fields}, {\bf 82}, 95--108.
 \item
Carne, T. K. (1985) A transmutation formula for Markov chains.
{\it Bull. Sci. Math.} {\bf 109}, 399--405.
 \item
Chandra, A.K., Raghavan, P., Ruzzo, W.I., Smolensky, R. and
Tiwari, P. (1989) The electrical resistance of a graph captures
its commute and cover times. {\it In Proc. 21 ACM Symp. Theory
of Computing,} 574--586.
 \item
Dembo A., Gantert N., Peres Y. \andd Zeitouni O. (2001) Large
deviations for random walks on Galton-Watson trees: Averaging
and Uncertainty. {\it Probab. Theory Relat. Fields}, To appear.
 \item
Lee, S. (1994a) Ph.D. thesis, Cornell University.
 \item
Lee, S. (1994b) Optimal drift on $[0,1]$. (English. English
summary) {\it Trans. Amer. Math. Soc.} {\bf 346}, 159--175.
 \item
Lyons, R. \andd Peres, Y. (2001) {\it Probability on Trees and
Networks} (a book), Cambridge University Press, in progress.
Current version published on the web at \\ {\tt
http://php.indiana.edu/\~{}rdlyons}.
 \item
Peres, Y. (1999) {\it Probability on Trees: An Introductory
Climb, Ecole d'Et\'e de Probabilit\'es de Saint Flour XXVII --
1997. Lecture Notes in Math 1717,} 193--280. Springer, Berlin.
 \item
Varopoulos, N. Th. (1985) Long range estimates for Markov
chains. {\it Bull. Sci. Math.} {\bf 109}, 225--252.
 \item
Vir\'ag, B. (2000) On the speed of random walks on graphs,
{\it Ann. Probab.} {\bf 28}, 379-394.
\end{enumerate}

\noindent 
Department of Mathematics  \\
Massachusetts Institute of Technology  \hfill {\tt balint@math.mit.edu} \\
Cambridge, MA 02139, USA  \hfill{\tt http://www-math.mit.edu/\~{}balint} 
\end{document}